\documentclass[12pt]{article}
\usepackage{amsmath,amsthm,amssymb}

\usepackage{epsf}
\usepackage{makeidx}

\theoremstyle{plain}
\newtheorem{theo}{Theorem}[section]
\newtheorem*{theo*}{Theorem}

\newtheorem*{lem*}{Lemma}

\theoremstyle{definition}
\newtheorem{remark}[theo]{Remark}
\newtheorem*{remark*}{Remark}

\newtheorem*{ex*}{Example}
\theoremstyle{remark}
\theoremstyle{definition}

\newtheorem*{defi*}{Definition}
\numberwithin{equation}{section}
\newcommand{\euc}[1]{\mathbb R^{#1}}
\newcommand{\wt}{\widetilde}

\newcommand{\Om}{\varOmega}
\newcommand{\imb}{\hookrightarrow}

\newcommand{\oW}{\mathaccent'27W}
\newcommand{\oH}{\mathaccent'27H}

\numberwithin{equation}{section}

\makeindex

\begin{document}
\title
{Linear elliptic differential equations}
\author{Ch\' erif Amrouche$^1$, Miroslav Krbec$^2$,\\ 
Brigitte Lucquin-Desreux$^3$, \v S\' arka Ne\v casov\' a$^2$}

\date{}
\maketitle

\renewcommand{\refname}{Further Reading}

{\parindent = 0pt $^1$Universit\' e de Pau et des Pays de l'Adour,
Laboratoire de Math\'ematiques Appliqu\' ees, I.P.R.A. - Av.\ de
l'Universit\' e, 64 000 Pau, France

$^2$Mathematical Institute, Academy of Sciences of the Czech Republic, 
\v Zitn\' a 25, 11567 Prague 1, Czech Republic

$^3$Laboratoire Jacques-Louis Lions, Universit\'e Pierre et
Marie Curie, 175 rue du Chevaleret, 75013 Paris, France
}

\bigskip

\noindent
{\bf Suggestions for cross-references to other articles within
Encyclopaedia:} 

{\it \parindent =0pt
251 - Partial differential equations in fluid mechanics,

158 - Viscous incompressible fluids: mathematical theory

430 - Inequalities in Sobolev spaces

388 - Functional analysis

318 - Classical mechanics
}

\section{Introduction}
\subsection{Motivation: A model problem}

Many physical problems can be modeled by partial
differential equations. Let us consider for example the case of an elastic 
membrane $\Om$, with fixed boundary $\Gamma$, subject to pressure forces $f$. 
The vertical membrane
displacement is represented by a real valued function $u$, which
solves the equation
\begin{equation}
\label{EqnLap}
-\Delta u (x) =  f (x), \qquad x =(x_1,x_2) \in
\Om,
\end{equation}
where the Laplace operator $\Delta$ is
defined, in two dimensions,  by
\begin{equation}
\Delta u = \frac{\partial^{2}u}{ \partial x^{2}_{1}} +
\frac{\partial^{2}u}{ \partial x_{2}^{2} }.
\nonumber
\end{equation}
As the membrane is glued to the curve $\Gamma $, $u$ satisfies the
condition
\begin{equation}
\label{ClDir}u(x)=0,  \quad x \in \Gamma.
\end{equation}
The system \eqref{EqnLap}-\eqref{ClDir} is the homogeneous
Dirichlet problem for the Laplace operator. It enters the more
general framework of (linear) elliptic boundary value problems, which
consist of a (linear) partial differential equation (in the example
above, of order two: the highest order in the derivatives) inside an
open set $\Om$ of the whole space $\euc N$, satisfying some
``elliptic'' property, completed by (linear) conditions on the
boundary $\Gamma$ of $\Om$, called ``boundary conditions''. In the
sequel, we only consider the linear case.

Our aim is to answer the following questions: does this problem
admit a solution? in which space? is this solution unique? does it
depend continuously on the given data $f$? In case of positive
answers, we say that the problem is ``well posed'' in the Hadamard
sense\index{well posed problem}. But other questions can also be
raised, such as the sign of 
the solution, for example, or its regularity.
We give a full survey of linear elliptic problems in a
bounded or in an exterior domain with a sufficiently
smooth boundary and in the whole space. In the general theory of
the elliptic problem we consider only smooth coefficients.
We survey the standard theory, which can be found in the
several well-known monographs of 1960s. The new trends in
investigation of the elliptic problems is to consider more
general domains with nonsmooth boundaries and nonsmooth
coefficients. On the other hand the regularity results for
elliptic systems have not been improved during last thirty
years. New trends also require employment of more general
function spaces and more general functional background. 

The number of references (Further reading) is strictly limited
here; there are only some of the most important
publications. The basic facts can usually be found in more
places and sometimes 
we do not mention the particular reference. To the very basic
references throughout belong 
\cite{Friedman}, \cite{GilbargTrudinger}, \cite{DautrayLions},
\cite{Ho}, \cite{Ladyzhenskaya}, \cite{LM}, \cite{RR}, \cite{W};
of course there are many others.



\subsection{The method}
To answer the above questions, we generally use, for such elliptic
problems, an approach based on what is called a ``variational
formulation'' (see section 3 below): the boundary value problem is
first transformed into a variational problem of lower order, which is
solved in a Hilbertian frame with help of the Lax-Milgram Theorem
(based on the representation theorem). All questions are then solved
(existence, uniqueness, continuity in terms of the data, regularity).
But this variational formalism does not necessary allow to treat all
the situations and it is limited to the Hilbertian case. Other
strategies can then be developed, based on a priori estimates and
duality arguments for the existence problem, or maximum principle for
the question of unicity. Without forgetting the particular cases
where an explicit Green kernel is computable (the Laplacean operator
in the whole space case for example).

Moreover, the study of linear elliptic equations is directly linked
to the function spaces background. It is the reason why we first
deal with Sobolev spaces---both of the integer and fractional order
and we survey their basic properties, imbedding and trace theorems.
We pay attention to the Riesz and Bessel potentials and we define
weighted Sobolev spaces important in the context of unbounded
opens. Second, we present the variational approach and the
Lax-Milgram theorem as a key point to solve a large class of boundary
value problems.  We give examples: the Dirichlet and Neumann problems
for the Poisson equation, the Newton problem for more general second
order operators; we also investigate mixed boundary conditions and
present an example of a problem of fourth order.  Then, we briefly
present the arguments for studying general elliptic problems and
concentrate on second order elliptic problems; we recall
the weak and strong maximum principle, formulate the Fredholm
alternative and tackle the regularity questions.  Moreover, we are
interested in the existence and uniqueness of solution of the Laplace
equation in the whole space and in exterior opens. Finally, we present
some particular examples arising from physical problems, either in
fluid mechanics (the  Stokes system) or in elasticity.

\section{Sobolev and other type of spaces}

Throughout, $\Om\subset \euc N$ will generally be an open
subset of the 
$N$-dimensional Euclidean space $\euc N$. 
A \textit{domain}\index{domain} 
will be an open and connected subset of $\euc N$.
We shall use standard
notations for the spaces $L^p(\Om)$, $C^{\infty}(\Om)$ 
etc.\ and their norms. Let us agree that $C^{k,r}(\Om)$, 
$k\in\mathbb N$, $r\in(0,1)$, denote the space of functions $f$ in
$C^{k}(\Om)$, whose derivatives $D^{\alpha} f$, 
$\alpha =(\alpha_1,\dots,\alpha_N) 
\in\mathbb N^N$ of order $|\alpha|=\sum_{i=1}^N \alpha_i=k$ are all 
$r$-H\"older continuous. 
By $\overline\Om$ in the notations for some of these
spaces we mean that the functions have the corresponding
property on $\Om$ and that they can be continuously extended to 
$\overline\Om$.

Let us recall several fundamental concepts. The space
$\mathcal{D}(\Om)$ of the \textit{test functions}\index{test
functions} in $\Om$ consists of all infinitely differentiable
$\varphi$ with a compact support in $\Om$. A locally convex topology
can be introduced here. The elements of the dual space
$\mathcal{D'}(\Om)$ are called the
\textit{distributions}\index{distribution}. If $f\in
L^1_{\textup{loc}}(\Om)$ (i.e.\  $f\in L^1(K)$ for all compact subset
$K$ of $\Om$), then $f$ is a \textit{regular distribution}; the
duality is represented by $\int_{\Om}f(x)\varphi(x)\,dx$.  If $f\in
\mathcal{D'}(\Om)$, we define the \textit{distributional} or the
\textit{weak derivative}\index{weak derivative} $D^{\alpha}$ of $f$
as the distribution $\varphi\mapsto (-1)^{|\alpha|}\langle
f,D^{\alpha}\varphi\rangle$.  Plainly, if $f\in L^1_{\textup{loc}}$
has ``classical'' partial derivatives in $L^1_{\textup{loc}}$, then
it coincides with the corresponding weak derivative.

If $\Om=\euc N$, it is sometimes more suitable to work with the
\textit{tempered distributions}. The role of $\mathcal{D}(\Om)$
is played by the space $\mathcal{S}(\euc N)$ of $C^{\infty}$
functions with finite pseudonorms $\sup |D^{\alpha}f(x)|(1+|x|)^{k}$,
$|\alpha|,k=0,1,2,\dots$. Recall that the Fourier transform
$\mathcal F$ maps $\mathcal{S}(\euc N)$ into itself and the same
is true for the space of the tempered distributions
$\mathcal{S}'(\euc N)$. 

\subsection{Sobolev spaces of positive order}

The \textit{Sobolev space}\index{Sobolev space!of integer order}
$W^{k,p}(\Om)$, $1\le p\le\infty$, 
$k\in\mathbb N$, is the space of all $f\in L^p(\Om)$, whose weak
derivatives up to order $k$ are regular distributions
belonging to $L^p(\Om)$; in $W^{k,p}(\Om)$ we introduce the norm
\begin{equation}\label{Sobnorm}
\|f\|_{W^{k,p}(\Om)}=\biggl(\sum_{|\alpha|\le k}\int_{\Om}
    |D^{\alpha}f(x)|^p\,dx\biggr)^{1/p}
\end{equation}
when $p<\infty$ and $\max_{|\alpha|\le k}\operatorname{sup\,ess}_{x\in\Om}
|D^{\alpha}f(x)|$ if $p=\infty$. 
The space $W^{k,p}(\Om)$ is a Banach space, separable for
$p<\infty$, reflexive for $1<p<\infty$; it is a Hilbert space
for $p=2$, more simply denoted $H^m(\Om)$. 
In the following we shall
consider only the range $p\in(1,\infty)$.

The link with the classical derivatives is given by this well-known
fact: A function $f$ belongs to $W^{1,p}(\Om)$ if and only if it
is a.e.\ equal to a function $\wt u$, absolutely continuous on almost
all line segments in $\Om$ parallel to the coordinate axes, whose
(classical) derivatives belong to $L^p(\Om)$ (the Beppo-Levi
theorem).

For $1<p<\infty$ and noninteger $s>0$ the Sobolev space
$W^{s,p}(\Om)$\index{fractional derivatives} of
order $s$ is\index{Sobolev space!of fractional order}
defined as the space of all $f$ with the finite norm
\begin{equation*}
\|f\|_{W^{s,p}(\Om)}= \biggl( \|f\|^p_{W^{[s]}_p(\Om)}
  +\sum_{|\alpha|=[s]}\int_{\Om}\!\int_{\Om}
    \frac{|D^{\alpha} f(x)-D^{\alpha}f(y)|^p}
     {|x-y|^{N+p(s-[s])}}\biggr)^{1/p},
\end{equation*}
where $[s]$ is the integer part of $s$.
Details e.g.\ in \cite{Ad}, \cite{Z}.

\subsection{Imbedding theorems}
One of the most useful and important features of the functions in
Sobolev spaces is an improvement of their integrability properties and
the compactness of various imbeddings. Theorems of this type were
first proved by Sobolev and Kondrashev\index{imbedding
theorems}\index{compact imbeddings}. Let us agree that
the symbol $\imb$ and $\imb\imb$ stands for an imbedding
and for a compact imbedding, respectively.

\begin{theo}
Let $\Om$ be a Lipschitz open. Then 
\begin{enumerate}
\item[\textup{(i)}] If $sp<N$, then $W^{s,p}(\Om)\imb L^{p^*}(\Om)$ with
$p^*=Np/(N-ps)$ (the Sobolev exponent). If $|\Om|<\infty$, then the
target space is any $L^r(\Om)$ with $0<r\le p^*$. 

If $\Om$ is
bounded, then $W^{s,p}(\Om)\imb\imb L^{q}(\Om)$ for all $1 \le q <p^*$.
\item[\textup{(ii)}]  
If $sp>N$, then $W^{j+s,p}(\Om)\imb C^{j}(\Om)$. If $\Om$
has the Lipschitz boundary, then $W^{j+s,p}(\Om) \imb
C^{j,\mu}(\overline\Om)$. 

If $sp>N$, then $W^{j+s,p}(\Om)\imb\imb
C^{j}(\Om)$ and $W^{j+s,p}(\Om)\imb\imb W^{j}_q(\Om)$ for all 
$1\le q\le \infty$. If, moreover, $\Om$ has the Lipschitz boundary,
then the target space can be replaced by $C^{j,\mu}(\overline \Om)$
provided $sp>N>(s-1)p$ and $0<\mu<s-N/p$.
\end{enumerate}
\end{theo}

Note that if the imbedding $W^{s,p}(\Om)\imb L^q(\Om)$ is compact for
some $q\ge p$, then $|\Om|<\infty$. Moreover, if
$\limsup_{r\to\infty}|\{x\in\Om;\, r\le |x|<r+1\}|>0$, then
$W^{s,p}(\Om)\imb L^q(\Om)$ cannot be compact.\hfill\break 

\subsection{Traces and Sobolev spaces of negative order}
Let $s>0$ and let $\Om$ be, for simplicity, a bounded open
subset of $\euc N$ with boundary $\Gamma$ of class $C^{[s],1}$.
Then with help of local coordinates, we can define  
Sobolev spaces $W^{s,p}(\Gamma)$ (also denoted $H^s(\Gamma)$ for
$p=2$) on $\Gamma = \partial\Om$ (see e.g.\ \cite{N},
\cite{Ad} for details). If $f\in C(\overline\Om)$,
then $f_{|\Gamma}$ has sense. Introducing the space ${\mathcal
D}({\overline \Om})$ of restrictions in $\Om$ of  functions in
${\mathcal D}(\euc N)$, one can show that if $f\in {\mathcal
D}({\overline \Om})$, we have
$\|f_{|\Gamma}\|_{W^{1-1/p,p}(\Gamma)}\le C\|f\|_{W^{1,p}(\Om)}$
so that, in view of the density of ${\mathcal D}({\overline
\Om})$ in $W^{1,p}(\Om)$, the restriction of $f$ to $\Gamma$ can
be uniquely extended to the whole $W^{1,p}(\Om)$. The result is
the bounded \textit{trace operator}\index{trace of a function}
$\operatorname{\gamma_0}: W^{1,p}(\Om)\to W^{1-1/p,p}(\Gamma)$.
Moreover, every $g\in {W^{1-1/p,p}(\Gamma)}$ can be extended to
a (non unique) function $f\in W^{1,p}(\Om)$ and this extension
operator is bounded with respect to the corresponding norms.

More generally, let us suppose $\Gamma$ is of class $C^{k-1,1}$ and
define the operator $\operatorname{Tr}_n$ for any $f \in {\mathcal
D}({\overline \Om})$ by 
${\operatorname{Tr}}_n f=(\gamma_0 f,\gamma_1f,\dots,
\gamma_{k-1} f),
$ where
\begin{equation*}
\gamma_j f(x) = \frac{\partial^jf}{\partial n^{j}}(x)=\sum_{|\alpha|=j}
 \frac{j !}{\alpha !}\,(\partial^{\alpha}f(x)/{\partial x^{\alpha}})
  n^{\alpha}, \qquad x\in\Gamma,
\end{equation*}
is the $j$-th order derivative of $f$ with respect to the outer
normal $n$ at $x\in\Gamma$; by density, this operator can be uniquely extended 
to a continuous linear mapping defined on the space $W^{k,p}(\Om)$; 
moreover, $\gamma_0 (W^{k,p}(\Om)) = {W^{k-1/p,p}(\Gamma)}$.

The kernel of this mapping is the space 
$\oW^{k,p}(\Om)$ (denoted by $H^k_0(\Om)$ for $p=2$), 
where $\oW^{s,p}(\Om)$ is 
defined as the closure of $\mathcal{D}(\Om)$ in
$W^{s,p}(\Om)$ ($s>0$). 
For $1<p<\infty$ the following holds: $\oW^{s,p}(\euc N)=W^{s,p}(\euc N)$,
$\oW^{s,p}(\Om)=W^{s,p}(\Om)$ provided $0<s\le 1/p$.
If $s<0$, then the space $W^{s,p}(\Om)$\index{negative
smoothness} is defined as the dual 
to $\oW^{-s,p'}(\Om)$, where $p'=p/(p-1)$ (see e.g.\
\cite{Tr1}, \cite{Tr2}). Observe that for an arbitrary $\Om$ a function $f\in
W^{1,p}(\Om)$ has the zero trace if and only if
$f(x)/\operatorname{dist}(x,\Gamma)$ belongs to $L^p(\Om)$.

For $p=2$, we  simply denote by $H^{-k}(\Om)$ the dual space of
$H^k_0 (\Om)$. In the case of  bounded opens, we recall the
following useful Poincar\'e-Friedrichs
inequality\index{Poincar\'e-Friedrichs  inequality} (for
simplicity, we state it here in the Hilbert frame):
 
\begin{theo}
Let $\Om$ be bounded (at least in one direction of the space). Then
there exists a positive constant $C_P(\Om)$ such that
\begin{equation}
\label{Eqn:Poincare}
\|v\|_{L_2(\Om)}\le C_P(\Om) \|\nabla v\|_{[L_2(\Om)]^N}
  \qquad \hbox{ for all } v \in H^1_0(\Om).
\end{equation}
\end{theo}
 
\subsection{The whole space case: Riesz and Bessel potentials}
The Riesz potentials $\mathcal I_{\alpha}$ naturally occur when
one defines the formal powers of the Laplace operator $\Delta$.
Namely, if $f\in\mathcal S(\euc N)$ and $\alpha>0$, then  
$\mathcal F\left[(-\Delta)^{\alpha/2}f\right](\xi)
=|\xi|^{\alpha}\mathcal Ff(\xi)$. This can be taken formally as
a definition of the \textit{Riesz potential} $\mathcal
I_{\alpha}$ on $\mathcal S'(\euc N)$,
\begin{equation*}
\mathcal I_{\alpha}f(.)
  =\mathcal F^{-1}\left[|\xi|^{-\alpha}\mathcal Ff(\xi)\right](.)
\end{equation*}
for any $\alpha\in\mathbb R$. If $0<\alpha<N$, then
$ I_{\alpha}f(x)
  =(I_{\alpha}*f)(x),
$
where $I_{\alpha}$ is the inverse Fourier transform of
$|\xi|^{-\alpha}$, 
\begin{equation*}
I_{\alpha}(x)=C_{\alpha}|x|^{\alpha-N}, \qquad
   C_{\alpha}=\Gamma\left((N-\alpha)/{2}\right)
    \left(\pi^{N/2}2^{\alpha}\Gamma(\alpha/2)\right)^{-1}
\end{equation*}
($\Gamma$ is the Gamma function) is the \textit{Riesz
kernel}.\index{Riesz kernel}
The following formula is also true:
\begin{equation*}
I_{\alpha}(x)=C_{\alpha}\int_{0}^{\infty}
   t^{\frac{\alpha-N}{2}}e^{-\frac{\pi |x|^2}{t}}
     \,\frac{dt}{t}.
\end{equation*}
Recall that every $f\in \mathcal S(\euc N)$ can be represented as
the \textit{Riesz potential} $\mathcal I_{\alpha}g$  of a
suitable function $g\in \mathcal S(\euc N)$, namely,
$g=(-\Delta)^{\alpha/2}f$; we get the representation formula
\begin{equation*}
f(x)=\mathcal I_{\alpha}g(x)
  =C_{\alpha}\int_{\euc N}\frac{g(y)}{|x-y|^{N-\alpha}}\,dy.
\end{equation*}
The standard density argument implies then an
appropriate statement for functions in $W^{k,p}(\euc N)$ with
an integer $k$ and for the \textit{Bessel potential
spaces} $\mathcal H^{\alpha,p}(\euc N)$---see below for their
definition. 
The original Sobolev imbedding theorem comes from the
combination of this representation and the basic continuity
property of $I_{\alpha}$, $\alpha p<N$, 
\begin{equation*}
I_{\alpha}: L^p(\euc N)\to L_q(\euc N), \qquad
           \frac{1}{q}=\frac{1}{p}-\frac{\alpha}{N}.
\end{equation*}

To get an isomorphic representation of a Bessel potential space 
(of a Sobolev space with positive integer smoothness in
particular) it is more convenient to consider the \textit{Bessel
potentials} (of order $\alpha\in\mathbb R$),\index{Bessel
potentials} 
\begin{equation*}
\mathcal G_{\alpha} f(x)=(G_{\alpha}*f)(x)
 =\mathcal F^{-1}\left([1+|\xi|^2]^{-
  \alpha/2}\mathcal Ff(\xi)\right)(x)
\end{equation*}
(with a slight abuse of the notations); the following formula for the
\textit{Bessel kernel}\index{Bessel kernel} $G_{\alpha}$ is
well known:
\begin{equation*}
G_{\alpha}(x)=c_{\alpha}^{-1}\int_{0}^{\infty}
   t^{\frac{\alpha-N}{2}}e^{-\frac{\pi |x|^2}{t}-\frac{t}{4\pi}}
     \,\frac{dt}{t}
\end{equation*}
(cf.\ the analogous formula for $I_{\alpha}$), where
$c_{\alpha}=(4\pi)^{\alpha/2}\Gamma(\alpha/2)$.
The kernels
$G_{\alpha}$ can alternatively be expressed with help of Bessel
or Macdonald functions.

Now we can define the \textit{Bessel potential spaces}. 
For $s\in\mathbb R$ and
$1<p<\infty$, let $\mathcal{H}^{s,p}(\euc N)$ be the space of all
$f\in\mathcal S'(\euc N)$ with the finite norm
\begin{equation*}
\|f\|_{\mathcal{H}^{s,p}(\euc N)}
 =\biggl(\int_{\euc N}\mathcal F^{-1}
    \left((1+|\xi|^2)^{s/2}\mathcal F f(\xi)\right)^p\,d\xi
      \biggr)^{1/p}.
\end{equation*}
In other words, the spaces $\mathcal{H}^{s,p}(\euc N)$ are
isomorphic copies of $L^p(\euc N)$.

For $k=0,1,2,\dots$, plainly $\mathcal{H}^{k,2}(\euc N)=W^{k,2}(\euc
N)$ by virtue of the Plancherel theorem. But it is true also for
integer $s$ and general $1<p<\infty$ (see e.g.\ \cite{Tr1}).

\begin{remark} 
Much more comprehensive theory of general Besov and
Lizorkin-Triebel spaces in $\euc N$ has been established in last
decades, relying on the the Littlewood-Paley theory. Spaces on
opens can be defined as restrictions of functions in the
corresponding space on the whole $\euc N$, allowing to
derive their properties from those valid for functions on $\euc
N$. The justification for that are extension theorems. In
particular, there exists a universal extension operator for the
Lipschitz open, working for all the spaces mentioned up to now.
We refer to 
\cite{Tr1}, \cite{Tr2}.
\end{remark}

\subsection{Unbounded opens and weighted spaces}\label{zerobellow}

The study of the elliptic problems in unbounded opens is usually
carried out with use of suitable Sobolev weighted space. 
The Poisson equation
\begin{equation}\label{eq5.1}
-\Delta u = f \quad\textup{in $\mathbb R^N$, $N \ge 2$},
\end{equation}
is the typical example; the Poincar\'e inequality
\eqref{Eqn:Poincare}  is not  
true here and it is suitable to introduce Sobolev spaces with
weights.  
 
Let $m\in \mathbb N$,
$1<p<\infty$, $\alpha\in \mathbb R$, $k= m-N/p -\alpha$ if $
{N}/{p}+ \alpha\in \{1,\dots,m\}$ and $k=-1$ elsewhere. For an
open $\Om \subset \mathbb R^N$ we define
\begin{equation*}
\begin{array}{rcl}
W^{m,p}_{\alpha}(\Om) &=& \{v\in{\mathcal D}'(\Om), \ 0
\le |\lambda| \le k, \ \rho^{\alpha-m-|\lambda |}(\log \rho)^{-1} D^{\lambda}
u\in L^p (\Om), \\
& \ & \ \ 
k+1 \le |\lambda| \le m, \ \rho^{\alpha-m+|\lambda|} D^{\lambda} u
\in L^p (\Om)\},
\end{array}
\end{equation*}
where $\rho(x) = (1+|x|^2)^{1/2}$, 
$\log\rho =\log (2+|x|^2)$. Note that
$W_{\alpha}^{m,p}$ is a reflexive Banach space for the norm 
$\|.\|_{W^{m,p}_{\alpha}}$ defined by
\begin{equation*}
\|u\|^p_{W^{m,p}_{\alpha}} =\sum_{0\leq |\lambda
 |\leq k} \|\rho ^{\alpha -m+|\alpha |} (\log\rho)^{-1}
  D^{\lambda }u\|^p_{L^p(\Om)} +\sum_{k+1\leq |\lambda
 |\leq m}\|\rho ^{\alpha -m+|\alpha |} 
  D^{\lambda }u\|^p_{L^p(\Om)}.
 \end{equation*}
We also introduce the following seminorm
\begin{equation*}
|u|_{W^{m,p}_{\alpha}} = \biggl( \sum _{|\lambda
 |=m }\|\rho ^{\alpha } D^{\lambda
  }u\|^p_{L^p(\Om)}\biggr)^{1/p} .
\end{equation*}

\noindent Let
\begin{equation*}
\oW_{\alpha}^{m,p} (\Om) = \{v\in W^{m,p}_{\alpha}; \,\gamma_0 (v)
= \dots \gamma_{m-1} (v) = 0\}.
\end{equation*}
If $\Om$ is a Lipschitz domain, then
$\oW_{\alpha}^{m,p} (\Om)$  is the closure of ${\mathcal D}(\Om) $ in 
 $W_{\alpha}^{m,p} (\Om)$, while 
 ${\mathcal D}(\overline{\Om})$ is dense in
$W_{\alpha}^{m,p}(\Om)$. We denote by
$W^{-m,p'}_{-\alpha}(\Om)$ the dual of $\oW_{\alpha}^{m,p}
(\Om)$. We note that these spaces contain also polynomials,
\begin{equation*}
P_j \subset W_{\alpha}^{m,p}(\Om) \Leftrightarrow
\begin{cases}
j= [m-\frac{N}{p}-\alpha] 
    & \textup{if $\frac{N}{p} + \alpha\notin Z$}\\ 
j=m-\frac{N}{p} -\alpha   &\textup{elsewhere},
\end{cases}
\end{equation*}
where $[s]$ is the integer part of $s$ and $P_{[s]}=\{0\}$
if $[s] < 0$. The fundamental property of functions belonging to
these spaces is that they satisfy the Poincar\'e weighted inequality.
A open $\Om $ is an \textit{exterior
domain}\index{exterior domain} if it is the complement of
a closure of a bounded domain in $R^N$. 

\begin{theo}\label{Th5.1}
Suppose that $\Om$ is an exterior domain or $\Om =
\mathbb R^N_+$ or $\Om = \mathbb R^N$. Then
\begin{enumerate}
\item[\textup{(i)}]
the seminorm $|\cdot|_{W^{m,p}_{\alpha}(\Om)}$
is a norms on $W^{m,p}_{\alpha} (\Om)/P_j$, equivalent to
the quotient norm with $j' = \min (m-1,j)$,
\item[\textup{(ii)}]
the seminorm $|\cdot|_{W^{m,p}_{\alpha}(\Om)}$
is equivalent to the full norm on $\oW^{m,p}_{\alpha}(\Om )$.
\end{enumerate}
\end{theo}

\section{Variational approach}
 Let us first describe the method on the model problem   
\eqref{EqnLap}-\eqref{ClDir}, supposing $f \in L^2(\Om)$ and
$\Om$ bounded. We first  
suppose that this problem admits a sufficiently smooth function
$u$. Let $v$ be any arbitrary (smooth) function; we multiply
equation \eqref{EqnLap} by $v(x)$ and integrate with respect to
$x$ over $\Om$; this gives 
\begin{equation*}
\int_\Om - (\Delta u v)(x)\,dx = \int_\Om (fv)(x) \,dx.
\end{equation*}
Using the following Green's formula\index{Green's formula}
($d\Gamma (x)$ denotes the measure on  
$\Gamma = \partial \Om$ and $\frac{\partial u}{\partial n}
(x) = \nabla u(x) \cdot n(x)$, where $n(x)$ is the unit normal at
point $x$ of $\Gamma$ oriented towards the exterior of $\Om$)
\begin{equation}
\label{Green}
\int_\Om(\Delta u v)(x)\,dx = - \int_\Om
(\nabla u \cdot \nabla v )(x) \,dx + \int_\Gamma
(\frac{\partial u}{\partial n} v )(\sigma ) \,d\sigma,
\end{equation}
we get, since $v|_\Gamma=0$: $ {\mathcal A}(u,v) = L(v)$, where we
have set 
\begin{equation}
\label{DefALLapDir}
{\mathcal A}(u,v)=\int_\Om \nabla u(x) \cdot \nabla v(x)\,dx, 
         \quad L(v)=\int_\Om f(x)v(x)\,dx .
\end{equation}
The idea is to study in fact this new problem (showing first its equivalence 
with the boundary value problem), noting
that it makes sense
for far less regular functions $u$, $v$ (and also $f$), in fact $u,v \in 
H^1_0(\Om)$ (and $f \in H^{-1}(\Om)$).

\subsection{The Lax-Milgram theorem}\index{Lax-Milgram theorem}

The general form of a variational problem\index{variational
problem} is 
\begin{equation}
\label{PV}
\mbox {to find } u \in V \ \mbox { such that }
\ {\mathcal A}(u,v) = L(v) \ \mbox{for all  } v \in V,
\end{equation}
where $V$ is a Hilbert space, ${\mathcal A}$ a bilinear continuous form
defined on $V \times V$ and $L$ a linear continuous form defined on
$V$.
We say moreover that ${\mathcal A}$ is
\textit{$V$-elliptic}\index{VB@$V$-ellipticity} if 
there exists a positive constant $\alpha$ such that
\begin{equation}
\label{ellipticite}
{\mathcal A}(u,u)
   \ge\alpha\| u \|_V^2  \quad \textup{for all $u \in V$}.
\end{equation}
The following theorem is due to Lax and Milgram.

\begin{theo}
Let $V$ be a Hilbert space. We suppose that ${\mathcal A}$ is a bilinear
continuous form on $V \times V$ which is $V$-elliptic and that $L$ is
a linear continuous form on $V$.  Then the variational problem
\eqref{PV} has a unique solution $u$ on $V$.  Moreover, if ${\mathcal A}$
is symmetric, $u$ is characterized as the minimum value on $V$ of the
quadratic functional $E$ defined by
\begin{equation}
\label{DefEnergie}
\hbox{ for all } v \in V, \quad E(v)
  =\frac{1}{2} {\mathcal A}(v,v)-L(v).
\end{equation}
\end{theo}

\begin{remark} 
\noindent (i) We have the following ``energy
estimate'':\index{energy estimate} 
$\|u\|_V \le \frac{1}{\alpha} \|L\|_{V'}$.
In the particular case of our model problem, this inequality shows
the continuity of the solution $u\in H^1_0(\Om)$ with respect to
the data $f \in L^2(\Om)$ (that can be weakened by choosing $f\in
H^{-1}(\Om)$).

\noindent (ii) Theorem 3.1 can be extended to sesquilinear continuous forms 
${\mathcal A}$ defined on 
$V \times V$; such form is called  $V$-elliptic if there exists a positive 
constant $\alpha$ such that
\begin{equation}
\label{ellipticiteC}
  {\rm Re} {\mathcal A}(u,u) 
   \ge \alpha \| u \|_V^2 \quad \textup{for all $u \in V$}.
\end{equation}
Finally, $V'$ is the dual of $V$.

\noindent (iii) Denoting by $A$ the linear operator defined on the
space $V$ by  ${\mathcal A} (u,v)=\langle Au,v\rangle_{V',V}$, for
all $v \in V$, the Lax-Milgram theorem shows that $A$ is an
isomorphism from $V$ onto its dual space $V'$, and the problem
\eqref{PV} is equivalent to solving the equation $Au=L$.
 
\noindent (iv) Let us do some remarks concerning the numerical
aspects. First, this variational formulation is the starting point of
the well known finite element method: The idea is to compute a
solution of an approximate variational problem stated on a finite
subspace of $V$ (leading to the resolution of a linear system), with
a precise control of the error with the exact solution $u$.
Second, the equivalence with a minimization problem allows the use
of other numerical algorithms. 
\end{remark} 

Let us now present some classical examples of second order elliptic
problems than can be solved with help of the variational theory.

\subsection{The Dirichlet problem for the Poisson equation}%
\index{Dirichlet problem} 
We consider the problem on a bounded Lipschitz open
$\Om\subset \euc N$ 
\begin{align}
\label{eq3.10}\begin{split}
-\Delta  u &=  f\  \qquad \textup{in $\Om$},\\
u &= u_0 \qquad \textup{on $\Gamma = \partial \Om$},
\end{split}
\end{align}
with $u_0 \in H^{1/2}(\Gamma)$, so that there exists  $U_0 \in H^1(\Om)$ 
satisfying $\gamma_0 (U_0) =u_0$.
The variational formulation of problem \eqref{eq3.10} is
\begin{equation}
\label{PVDirNH}
\hbox {to find } u \in U_0 + H^1_0(\Om) \hbox { such that for all  } v\in 
H^1_0(\Om),
\quad {\mathcal A}(u,v)=L(v),
\end{equation}
with ${\mathcal A}$ given by \eqref{DefALLapDir} and a more general $L$ 
with $f \in H^{-1}(\Om)$, defined by
\begin{equation}
\label{Eqn:DefnLGene}
L(v)=\langle f,v\rangle _{H^{-1}(\Om),H^1_0(\Om)}.
\end{equation}
The existence and uniqueness of a solution of \eqref{PVDirNH}
follows from Theorem 3.6 (and Poincar\'e inequality
\eqref{Eqn:Poincare}. Conversely, thanks to the density of ${\mathcal
D}(\Om)$ in $H^1_0(\Om)$, we can show that $u$ satisfies
\eqref{eq3.10}. More precisely, we get:

\begin{theo}
\label{thmNHDirLapcu}
Let us suppose $f \in H^{-1}(\Om)$ and $u_0 \in H^{1/2}(\Gamma)$; let $U_0
\in H^1(\Om)$
satisfy $\gamma_0(U_0) =u_0$. Then the boundary value problem
\eqref{eq3.10} has a unique solution $u$
such that
$u-U_0 \in H^1_0(\Om)$. This is also the unique solution of the
variational problem \eqref{PVDirNH}. Moreover, 
there exists a positive constant $C= C(\Om)$ such that
\begin{equation}
\| u \|_{H^1 (\Om)} 
 \le C \left(\| f \|_{H^{-1}(\Om)}
 +\| u_0 \|_{H^{1/2}(\Gamma)} \right),
\end{equation}
which shows that $u$ depends continuously on the data $f$ and $u_0$.
\end{theo}

Moreover, using technics of Nirenberg's differential quotients, 
we have the following regularity result (see e.g.\ \cite{Grisvard}):%
\index{regularity}

\begin{theo}
Let us suppose that $\Om$ is a bounded open subset of $\euc N$ with a
boundary
of class ${\mathcal C}^{1,1}$ and let satisfy $f \in L^2(\Om)$,
$u_0 \in H^{3/2}(\Gamma)$.
Then $u \in H^2(\Om)$ and each equation in \eqref{eq3.10}
 is satisfied almost
everywhere (on $\Om$ for the first one and on $\Gamma$ for the boundary 
condition). Moreover, there exists a
positive constant
$C = C(\Om)$ such that
\begin{equation}
\| u \|_{H^2(\Om)} \le C \, [ \| f
\|_{L^2(\Om)}
+\| g \|_{H^{3/2}(\Gamma)}].
\end{equation}   
\end{theo}

By induction, if the data are more regular, i.e.\  $f \in H^k
(\Om)$ and $u_0 \in H^{k + 3/2}(\Gamma)$ (with $k \in  \mathbb
N$), and  if $\Gamma$ is of class  ${\mathcal C}^{k+1,1}$, we get $u \in
H^{k+2} (\Om)$.

\begin{remark} 
Let us point out the importance of the open geometry.  For example,
if $\Om$ is a bounded plane polygon, one can find $u \in
H^1_0(\Om)$ with $\Delta u \in {\mathcal C}^{\infty} ({\overline \Om})$,
such that $u \notin H^{1+\pi/w}(\Om)$, where $w$ is the biggest
value of the interior angles of the polygon.  In particular, if the
polygon is not convex, the solution of the Dirichlet problem
\eqref{eq3.10} cannot be in $H^2(\Om)$.
\end{remark} 

\subsection {The Neumann problem for the Poisson equation}
\index{Neumann problem}
We consider the problem ($n$ is the unit outer normal on $\Gamma$)
\begin{align}
\label{eq3.11}
\begin{split}
-\Delta  u &= f \qquad \textup{in $\Om$} \\
\frac {\partial u}{\partial n}&= h
         \qquad \textup{on $\Gamma$}.
\end{split}
\end{align}
Setting $E(\Delta) =\{v\in H^1(\Om)$;\,  $\Delta v\in L^2(\Om)\}$,
the space ${\mathcal D}({\overline \Om})$ is a dense subspace, and we
have the following Green formula for all $u \in E(\Delta)$ and $v \in
H^1(\Om)$:
\begin{equation*} 
\int_{\Om}  \Delta u (x) v(x) \, dx=-\int_{\Om} \nabla u (x) 
\cdot \nabla v (x) \,dx 
 +\langle {\frac{\partial u}{\partial n}},\gamma_0 
  v \rangle_{H^{-1/2}(\Gamma), H^{1/2}(\Gamma)}.
\end{equation*}
If $u \in H^1(\Om)$ satisfies \eqref{eq3.11} with $ f \in
L^2(\Om)$ and $h
\in H^{-1/2}(\Gamma)$, then for any function $v \in H^1(\Om)$, we have, 
by virtue of the above Green formula,
\begin{equation*}
{\mathcal A}(u,v) = {\tilde L}(v), \quad  {\tilde L} v=\int_{\Om} (f v) 
(x)\,dx+\langle h,\gamma_0 v \rangle _{H^{-1/2}(\Gamma), H^{1/2}(\Gamma)}.
\end{equation*}
But here, the form ${\mathcal A}$ is not $H^1(\Om)$-elliptic; in fact,
one can check that, if problem \eqref{eq3.11} has a solution, then we
have necessarily (take $v=1$ above)
\begin{equation}
\label{Eqn:CompDonWeak}
\int_{\Om}f(x)\,dx+\langle h,1\rangle_{H^{-1/2}(\Gamma),H^{1/2}(\Gamma)}=0.
\end{equation}
Moreover, we note that if $u$ is a solution, then $u +C$, where $C$
is an arbitrary constant, is also a solution. So the variational
problem is not well posed on $H^1(\Om)$.  It can be, however, solved
in the quotient space\index{quotient space} $H^1(\Om)/{\mathbb
R}$ which is a Hilbert space for the quotient
norm\index{quotient norm} 
\begin{equation}\label{eqnormequotient}
\| \dot v \|_{H^1(\Om)/{\mathbb R}} 
   = \inf_{k \in \mathbb R} \| v + k \|_{1,\Om},
\end{equation}
but also for the semi-norm $v\mapsto |v|_{H^1(\Om)}=
{\sqrt{{\mathcal A}(v,v)}}$, which is an equivalent norm on this
quotient space, see~\cite{N}.

\noindent
Then, supposing that the data $f$ and $h$ satisfy the ``compatibility
condition'' \eqref{Eqn:CompDonWeak}, we can apply the Lax-Milgram
theorem to the variational problem
\begin{equation}\label{eq3.12}
\textit{to find $\dot u \in
H$ such that ${\mathcal A}(\dot u,\dot v) 
            = {\tilde L} (\dot v)$ for all $\dot v \in V$}
\end{equation} 
with $V=H^1(\Om)/{\mathbb R}$.
We get the following result 
(see e.g.\ \cite{N}):

\begin{theo}
Let us suppose that $\Om$ is
connected and that the data $f \in L^2(\Om)$ and $h \in
H^{-1/2}(\Gamma)$  satisfy \eqref{Eqn:CompDonWeak}.
Then the variational problem
\eqref{eq3.12} has a unique solution $\dot u$ in the space
$H^1(\Om)/{\mathbb R}$ and this solution is continuous with
respect to the data, i.e.\ there exists a positive
constant $C = C(\Om)$ such
that
\begin{equation*}
 | u |_{H^1(\Om)} \le  C 
\left(\|f\|_{L^2(\Om)}  + \|h \|_{H^{-1/2}(\Gamma)}\right)
    \quad \textup{for all $u \in \dot u$}.
\end{equation*}
\end{theo}
Moreover, if $\Gamma$ is of class  ${\mathcal C}^{1,1}$ and if  the
data satisfy $f \in L^2(\Om)$,
$g \in H^{1/2}(\Gamma)$, then every $u \in \dot u$ is such that $u \in 
H^2(\Om)$ and it satisfies each equation in \ref{eq3.11} almost 
everywhere.\index{regularity}

\subsection{Problem with mixed boundary conditions}%
\index{mixed boundary conditions}
Here we consider more general boundary conditions: the Dirichlet
conditions on a closed subset $\Gamma _{1}$ of $\Gamma =\partial
\Om$, and the Neumann, or more generally the ``Robin'',
conditions on the other part $\Gamma _2 = \Gamma -\Gamma _{1}$. We
seek $u$ such that ($f \in L^2(\Om)$, $h \in L^2(\Gamma_2)$, $a
\in L^{\infty}(\Gamma_2)$)
\begin{align}\label{3.15}\begin{split}
-\Delta  u &= f \qquad \ \textup{in $\Om$,} \\
u&= 0 \qquad \ \textup{on $\Gamma _1$}, 
 \quad a u + \frac{\partial u}{\partial n}= h 
   \qquad \textup{on $\Gamma_2$}.
\end{split}
\end{align}
Let $ V=\{v \in H^1(\Om );\, \gamma_0 v=0 \mbox { on } \Gamma _1\}$.
Then \eqref{PV} is the variational formulation of this problem with
\begin{enumerate}
\item[\textup{(i)}]
${\mathcal A}(u,v) = \int_\Om \nabla u(x) \cdot \nabla v(x) \,dx
  +\int_{\Gamma_2} (a \gamma_0 u \, \gamma_0 v)(\sigma) \,d\sigma,$
\item[\textup{(ii)}]
$L(v)=\int_\Om f(x) v(x)\,dx + \int_{\Gamma_2} (h\gamma_0 v)(\sigma)
\, d\sigma$.
\end{enumerate} 
Supposing for example $a \ge 0 $, we get a unique solution $u \in V$
for this variational problem by virtue of the Lax-Milgram theorem.
Moreover, if $u \in H^2(\Om)$, then $u$ is the unique solution in
$H^2(\Om) \cap V$ of the problem \eqref{3.15}.

\subsection { The Newton problem for more general operators}%
\index{Newton problem}
Let $\Om$ be a bounded open subset of $\euc n$.
We now consider more general second order operators of the form 
$v\mapsto -\nabla . (M \nabla v)+ b \cdot \nabla v+ c v,
$ where $b \in [W^{1,\infty}(\Om)]^N$, $c \in L^{\infty}(\Om)$,
$M$ is an $N \times N$ square matrix with entries $M_{ij}$, and
$\nabla \cdot (M \nabla v)$ stands for 
$\sum_{i,j =1}^N  \frac{\partial}{\partial x_{i}}[M_{ij}
\frac{\partial u}{\partial x_{j}}]$.
We also assume that there is a positive  constant $\alpha_M$ such
that
\begin{equation*}
\sum_{i,j=1}^N M_{ij}(x) \xi_i \xi_j \ge \alpha_M
\sum_{i=1}^N \xi_i^2
\hbox{ for a.e.\ } x \in \Om  \hbox{ and } \xi=(\xi_1,...,\xi_N) \in
\euc N.
\end{equation*}
For given data $f \in L^2(\Om), h \in L^2(\Gamma)$, we look for
$u$ solution of the problem
\begin{align}\label{eq3.14}\begin{split}
-\nabla \cdot (M\nabla u) + b \cdot \nabla u + cu &=  f \qquad \hbox{ in } 
\Om, \\
au+n\cdot(M\nabla u)&=h \quad \hbox{ on } \Gamma.
\end{split}
\end{align}
We assume that $a \in L^{\infty }(\Gamma)$. The variational
formulation of this problem is still \eqref{PV},  with $V =
H^1(\Om)$ and 
\begin{equation}
\label{DefAGenSymMixcu}
{\mathcal A}(u,v) = \int_\Om  M \nabla u  \cdot  \nabla v \,dx + 
\int_\Om [b \cdot \nabla u + c u]  v \,dx +  \int_{\Gamma} a \gamma_0 u 
\gamma_0 v \,d\sigma,
\end{equation}
 \begin{equation}
\label{DefLGenSymMixcu}
\quad L(v)=\int_\Om f(x) v(x) \,dx+\int_{\Gamma}(h\gamma_0 v)(\sigma)
\, d\sigma,
\end{equation}
If the following conditions 
\begin{equation*}
c -\frac{1}{2} \nabla \cdot b \ge C_0 \ge 0~ \ \hbox{ a.e.\ on } 
\Om, \quad 
a+\frac{1}{2} b \cdot \nu ge C_1 \, \ge 0~ \ \hbox{ a.e.\ on } 
\Gamma
\end{equation*}
are fulfilled, with $(C_0,C_1) \ne (0,0)$, then the bilinear form
${\mathcal A}$ is V-elliptic and the Lax-Milgram theorem applies.

\subsection{A biharmonic problem}\index{biharmonic problem}

We consider the Dirichlet problem for the operator of fourth order:
 ($c \in L^{\infty}(\Om)$):
  \begin{equation}
  \label{EqnBiLapcu}
\Delta^2 u + c u=f \quad \hbox{ in } \Om,
  \end{equation}
  \begin{equation}
  \label{ClDirNHBiLap}
 u= u_0 \ \hbox{ on } \Gamma, \quad 
 {\frac{\partial u}{\partial n}}= h  \ \hbox{ on }  \Gamma.
  \end{equation}
\begin{theo}
\label{thmNHDirBiLapcu}
Let us suppose that $\Om$ has a boundary
of class ${\mathcal C}^{1,1}$ and that the data satisfy
 $f \in H^{-2}(\Om)$, $u_0 \in H^{3/2}(\Gamma)$, $h \in H^{1/2}(\Gamma)$.
Let $U_0
\in H^2(\Om)$ be such that $\gamma_0(U_0) =u_0, \gamma_1(U_0)
=h$. Then, if  $c \ge 0$ a.e.\ in $\Om$, the boundary
value problem \eqref{EqnBiLapcu}-\eqref{ClDirNHBiLap} has a unique
solution $u$ such that
$u-U_0 \in H^2_0(\Om)$, and $u$ is also the unique solution of the
variational problem
\begin{equation}
\label{PVBiLapcuDirNH}
\hbox {to find } u\in U_0 + H^2_0(\Om)\ \hbox {such that\ } 
{\mathcal A}(u,v) =l(v) \ \mbox{for all\ } v\in H^2_0(\Om),
\end{equation}
where $l(v) = \langle f,v \rangle_{H^{-2}(\Om), H^{2}_0(\Om)}$ and
\begin{equation}
\label{DefALBiLapDircuu}
{\mathcal A}(u,v) = \int_\Om \Delta u(x)  \Delta v(x) \,dx
  +\int_\Om (c u v)(x) \,dx.
\end{equation}
Moreover, there exists a positive constant $C=C(\Om)$ such that
\begin{equation}
\| u \|_{H^2(\Om)}\le C \, [\| f
\|_{H^{-2}(\Om)}
 +\| u_0 \|_{H^{3/2}(\Gamma)}
  +\| h \|_{H^{1/2}(\Gamma)}],
\end{equation}
which shows that $u$ depends continuously upon the data $f$, $u_0$ and $h$.
\end{theo}

\begin{remark}
The Hilbert space choice $V$ is of crucial importance for the 
$V$-ellipticity. In fact, let us consider for example the  problem
\eqref{EqnBiLapcu}, \eqref{ClDirNHBiLapMod}, with
  \begin{equation}
  \label{ClDirNHBiLapMod}
 \Delta u= 0 \ \hbox{ on } \Gamma, \quad 
 {\frac{\partial\Delta u}{\partial n}}= 0  \ \hbox{ on }  \Gamma.
  \end{equation}
In fact, the associated bilinear form is not $V$-elliptic for 
$V=H^2(\Om)$ but it is $V$-elliptic for 
$V=\{v \in L^2(\Om);\, \Delta v \in L^2\Om)\}$.
\end{remark}

\section{General elliptic problems}\index{general elliptic problems}

Here  $\Om$ will be a bounded and sufficiently regular open subset
of $\euc N$. Let us consider a general linear differential operator
of the form
\begin{equation}\label{Def:OpGeneNC}
  A(x, D) u=\sum_{|\mu| \le l} a_{\mu}(x) D^{\mu} u, \quad a_{\mu}(x) 
   \in \mathbb C
\end{equation}
Setting $A_0 (x, \xi) = \sum_{|\mu| = l} a_{\mu}(x)
\xi^{\mu} $, we say that the operator $A$ is \textit{
elliptic}\index{operator!elliptic} at a point
$x$ if  $ A_0  (x, \xi) \ne 0$ for all $\xi \in \euc N -\{0\}$. One
can show that, if $N \ge 3$, $l$ is even, i.e.\  $l=2m$; the same
result holds for $N=2$ if the coefficients $a_{\mu}$ are real.
Moreover, for $N \ge 3$, every elliptic operator is \textit{properly
elliptic},\index{operator!properly elliptic} in the following sense: For
any independent vectors $\xi$, 
$\xi'$ in $\euc N$, the polynomial $\tau \mapsto  A_0 (., \xi + \tau
\xi')$ has $m$ roots with positive imaginary part.
 
The aim here is to study boundary value problems of the following type:
\begin{equation}\label{EqnGeneElliptic}
Au=f  \qquad \hbox{ in } \Om,
\end{equation}
\begin{equation}\label{ClGeneElliptic}
B_j u= g_j \, \mbox{ on } \Gamma, \quad j = 0,...,m-1,
\end{equation}
where $A$ is properly elliptic on ${\overline \Om}$, with sufficiently
regular coefficients, and the operators $B_j$ are boundary operators,
of order $m_j \le 2m-1$, that must satisfy some compatibility
conditions with respect to the operator $A$ (see
\cite{RR} for details; these conditions were introduced by Agmon,
Douglas and Nirenberg).  For example, $A = (-1)^m \Delta ^m$ and $B_j
= {\frac{\partial ^j}{\partial n^j}}$ is a convenient choice.

In order to show that problem
\eqref{EqnGeneElliptic}-\eqref{ClGeneElliptic} has a solution $u
\in H^{2m+r}(\Om)$ ($ r \in \mathbb N$), the idea is to show
that the operator ${\mathcal P}$ defined by $u \mapsto {\mathcal P} (u)=
(Au, B_0u,..., B_{m-1} u) $ is an index operator from
$H^{2m+r}(\Om)$ into $G = H^{r}(\Om) \times
\Pi_{j=0}^{m-1} H^{2m+r-m_j-1/2}(\Gamma)$  and to express the
compatibility conditions through the adjoint problem.

We recall that a linear continuous operator ${\mathcal P}$ is an
index operator if

\noindent{(i)}  dim Ker ${\mathcal P} <  \infty$, and Im ${\mathcal P}$ closed

\noindent{(ii)}  codim Im ${\mathcal P} < \infty$.

Then the index $\chi ({\mathcal P})$ is given by $\chi ({\mathcal P}) = $ dim 
Ker 
${\mathcal P}$ - codim Im ${\mathcal P}$.
We recall the following Peetre's theorem:

\begin{theo}
Let $E$, $F$ and $G$ be three reflexive Banach spaces such that $E
\imb\imb F$, and ${\mathcal P}$ a linear continuous operator from $E$ to
$ G$. Then condition (i) is equivalent  to

\noindent{(iii)} there exists $C \ge 0$, such that for all $u \in E$, we have
\begin{equation*} 
\|u\|_{E}\le C\,\left(\|{\mathcal P} u \|_{G}+\|u\|_{F} .
\right)
\end{equation*}
\end{theo}

Applying this theorem to our problem 
\eqref{EqnGeneElliptic}-\eqref{ClGeneElliptic}, condition (i) results
from a  priori estimates of the following type:
\begin{equation*} 
\|u\|_{H^{2m+r}(\Om)} 
  \le C \, \left (\|{\mathcal P} u  \|_{G}  
  +  \| u \|_{H^{2m+r-1}(\Om)}\right )
\end{equation*}
and condition (ii) by similar a priori estimates for the dual problem.

\section{Second order elliptic problems}
We consider a second order differential operator of the
``divergence form''\index{operator!elliptic of divergence form} 
\begin{equation}\label{eq3.4}
A u = -\sum^N_{i,j=1} (a^{ij}(x) u_{x_i})_{ x_j} + \sum^N_{i=1} b^i
(x) u_{x_i} + c(x) u 
\end{equation}
with given coefficient functions $a^{ij},b^i,c$ $(i,j=1,\dots,N)$, and
where we have used the notation $u_{x_i} = {\frac{\partial
u}{\partial x_i}}$.  
Such operators are said \textit{uniformly strongly elliptic}
\index{operator!uniformly strongly elliptic} in $\Om $
if there exists $\alpha >0$ such that
\begin{equation*} 
\sum_{|i|=|j|= 1}a^{ij}(x)\xi^{i}\xi ^j \geq \alpha |\xi|^{2}
\quad   \hbox{for all }   x\in\Om,\ \xi\in\euc N.
\end{equation*}

\begin{remark}
There exists elliptic problems for which the associated variational
problem does not necessarily satisfy the ellipticity condition. Let
us consider the following example, due to Seeley: Let $\Om =
\{(r,\theta) \in (\pi, 2 \pi) \times [0,2 \pi] \}$ and $A = - (e^{i
\theta} {\frac{\partial}{\partial \theta}})^2 - e^{2 i \theta} (1 +
{\frac{\partial^2}{\partial r^2}})$. One can check that, for all
$\lambda \in \mathbb C$, the problem $Au+ \lambda u = f$ in $\Om$
and $u=0$ on $\Gamma$ admits nonzero solutions $u$ which are given
by (with $\mu $ such that $\mu^2 = \lambda$) $u = \sin r \cos ( \mu
e^{- i \theta}) $ and $u = \sin r sin ( \mu e^{- i \theta}) $ for
$\lambda \ne 0$; $u = \sin r$ and $u = \sin \theta  e^{- i \theta} $
for $\lambda =0$.
\end{remark}

Most of results concerning existence, unicity, regularity for second
order elliptic problems can be established thanks to a maximum
principle. There exist different types of maximum principles, that we
now present.

\subsection{Maximum principle}

\begin{theo}[weak maximum principle]\label{Th2.121}
\index{weak maximum principle}
Let $A$ be a uniformly strongly elliptic operator of the form
\eqref{eq3.4} in a bounded open $\Om \subset
\euc N$, with $a^{ij}, b^i, c
\in L^{\infty} {(\Om})$ and $c \ge 0$.
 Let $u\in C^2(\Om) \cap C (\overline {\Om})$ and  
\begin{equation*}
A u \ge 0 \ [\hbox{resp. } A u \leq 0] \ \mbox{in} \ \Om.
\end{equation*}
Then
\begin{equation*}
\inf_{\Om} u  \ge 
\inf_{\partial\Om} u^{-} \ [\hbox{resp. } 
   \sup_{\Om} u \le  \sup_{\partial\Om} u^{+}],
\end{equation*}
where $u^{+}=\max(u,0)$ and $u^{-}= -\min(u,0)$. 
If  $c = 0$ in $\Om$, one can replace $u^{-}$\  $[$resp.\ $u^{+}]$ by
$u$. 
\end{theo}

\begin{theo}[strong principle maximum]\label{Th2.131}
\index{strong principle maximum} Under the assumptions of the above
theorem, if $u$ is not a constant function in $ C^2(\Om) \cap C
(\overline {\Om})$ such that $Au \ge 0$  [{ resp. } $ Au \leq 0
$] , then $\inf_{\Om} u  < u(x)$ [resp. $\sup_{\Om} u >
u(x) $], for all $x
\in \Om$.
\end{theo}

\begin{remark}
These two maximum principles can be adapted to elliptic
operators in nondivergence form,\index{operator!elliptic of
nondivergence form} i.e., 
\begin{equation}\label{eq3.5}
A u = -\sum^N_{i,j=1} a^{ij} (x) u_{x_i x_j} + \sum^N_{i=1} b^i
(x) u_{x_i} + c(x) u.
\end{equation}

\end{remark}

\subsection{Fredholm alternative}\index{Fredholm alternative}

We now present some existence results which are based on on the
Fredholm alternative rather than on the variational method.

Let us consider two Hilbert spaces $V$ and $H$, where $V$ is a dense
subspace of $H$ and $V \imb \imb H$.  Denoting by $V'$ the dual
space of $V$, and identifying $H$ with its dual space, we have
the following imbeddings: $V \imb H \imb V'$. Let ${\mathcal A}$ be a
sesquilinear form on $V \times V$,
$V$-coercive\index{coercivity} with respect to $H$, 
that is, there exists $\lambda_0 \in \mathbb R$ and $\alpha >0$ such that
\begin{equation*} 
 {\rm Re}  ( {\mathcal A} 
(v,v)) + \lambda_0 \|v\|^2_H
 \ge \alpha \|v\|^2_V \quad \textup{for all $v \in V$}.
\end{equation*}

Denoting by $A$ the operator associated with the bilinear form ${\mathcal
A}$ (see item (iii) of Remark 3.1), the equation $Au=f$ is equivalent
to $u-\lambda_0Tu =g$, with $ T = (A + \lambda_0 Id)^{-1} $ and $g =
T f$. Note that $T$ is an isomorphism from $H$ onto
$D(A)=\{u\in H;\, Au \in H \}$).
 
The operator $T: H \to H$ is compact and, thanks to the Fredholm
alternative, there are two situations:

\noindent{(i)} either Ker ${A} =0$ and ${A}$ is an isomorphism from
$D({ A})$ onto $H$

\noindent{(ii)} or Ker ${ A} \ne 0$; then  Ker $ {A}$ is of finite
dimension, and the problem $Au=f$ with $f \in H$ admits a solution if
and only if $f \in {\rm Im} { A} = [ {\rm Ker} ( {A^*}) ]^{\perp}$.

We now give another example in a non Hilbertian frame.  Let us
consider the problem \cite{Grisvard}: $Au=f$ in $\Om$ and 
$Bu=g $ on $\Gamma$ where $ \Gamma$ is of class $ C^{1,1}$, $A$,
which is defined by \eqref{eq3.4}, is uniformly strongly elliptic
with $ a^{ij} = a^{ji } \in C^{0, 1} (\overline \Om)$, $ b^i, c
\in L^{\infty}(\Om)$ and $ Bu = \gamma_0(u)$ or $ Bu =
\gamma_1(u)$. One can show that the operator $u\mapsto (A u, Bu)$ is a
Fredholm operator of index zero from $ W^{2, p}(\Om)$ in $
L^p(\Om) \times W^{2 - d - 1/p, p}(\Gamma)$  (with $d = 0$ if $ Bu
= \gamma_0(u)$ and $ d = 1$ if $Bu = \gamma_1(u)$).

\subsection{Regularity}
\index{regularity}

Assume that $\Om $ is a bounded open. Suppose that
$u\in H^1_0(\Om)$ is a weak solution of the equation
\begin{align}
\label{eq3.10Gene}
\begin{split}
A u & = f  \qquad \textup{in $\Om$},\\
u &= 0     \qquad \textup{on $\Gamma$}, 
\end{split}
\end{align}
where $A$ has the divergence form \eqref{eq3.4}.
We now address the question whether $u$ is in fact smooth: this is the
\textit{regularity problem} for weak solutions.

\begin{theo}[$H^2$-regularity] Let $\Om$ be open, of class $C^{1,1}$,
 $a^{ij}\in C^1(\overline{\Om})$, 
$b^i, c\in L^{\infty}(\Om )$, $f\in L^2(\Om)$.
Suppose furthermore that $u\in H^1(\Om)$ is a weak
solution of \eqref{eq3.10Gene}. Then 
$u\in H^2(\Om)$ and 
we have the estimate
\begin{equation*}
\|u\|_{H^2(\Om)}\leq C(\|f\|_{L^2(\Om)}
    +\|u\|_{L^2(\Om)}),
\end{equation*}
where the constant $C$ depends only on $\Om$ and on the
coefficients of $A$.
\end{theo}

\begin{theo}[higher regularity]
Let m be a nonnegative integer, $\Om$ be open, of class $C^{m + 1,1}$ 
and assume that
$a^{ij}\in C^{m+1}(\overline \Om)$, $b^i, c\in
C^{m+1}(\overline \Om)$, $f\in H^{m}(\Om)$.
Suppose furthermore that $u\in H^1(\Om)$ is a weak
solution of \eqref{eq3.10Gene}. Then 
$u\in H^{m+2}(\Om)$ and
\begin{equation*}
\|u\|_{H^{m+2}(\Om)}\leq C ( \|f\|_{H^{m}(\Om)}+\|u\|_{L^2(\Om)}),
\end{equation*}
where the constant $C$ depends only on $\Om$ and on the
coefficients of $A$.
In particular, if $m>{N}/{2}$, then $u\in C^2(\overline \Om
)$. Moreover if $\Om$ is of $C^{\infty }$ class and $f\in
C^{\infty}(\overline {\Om})$,  $a^{ij}\in 
C^{\infty}(\overline{\Om})$, 
$b^i, c\in C^{\infty}(\overline{\Om} )$, then $u \in
C^{\infty}(\overline {\Om})$.
\end{theo}

\begin{remark}

\noindent{(i)} If $u \in H^1_0(\Om)$ is the unique solution of 
\eqref{eq3.10Gene}, one can omit the $L^2$ norm of $u$ in the
right hand side of the above estimate. 
 
 \noindent{(ii)} Moreover, let us suppose the coefficients $a^{ij},
b^{i}$ and $ c$ are all $C^{\infty}$ and $f \in C^{\infty}(\Om)$;
then, if $u \in H^1(\Om)$ satisfies $Au =f$, $u \in
C^{\infty}(\Om)$; this is due to the ``hypoellipticity'' property
satisfied by  the operator $A$.  
\end{remark}

We have a similar result in the $L^p$ frame \cite{Grisvard}:

\begin{theo}[$W^{2,p}$-regularity] Let $\Om$ be open, of class $C^{1,1}$,
 $a^{ij}\in C^1(\overline{\Om})$, 
$b^i, c\in L^{\infty}(\Om )$.
Suppose furthermore that $b^i = 0, 1\leq i\leq N $ and $c \geq 0$ a.e or $ c 
\geq \beta > 0$  a.e . 
Then for every $f\in L^p(\Om)$ there exists a unique solution
$u\in W^{2,p}(\Om)$ of \eqref{eq3.10Gene}.
\end{theo}

\section{Unbounded open}

\subsection{The whole space}

Note in passing that we shall work with the weighted
Sobolev spaces $W^{m,p}_{\alpha}(\Om)$ defined in
Subsection~\ref{zerobellow}.  

\begin{theo}\label{Th5.2} The following claims hold true:
\begin{enumerate}
\item[\textup{(i)}]
Let $f\in W^{-1,p}_0(\mathbb R^N)$ satisfy the
compatibility condition
\begin{equation*}
\langle f,1\rangle_{W^{-1,p}_0(\mathbb R^N)\times W^{1,p'}_0(\mathbb
R^N)} = 0 \quad\mbox{if}\quad p' \ge N.
\end{equation*}
Then the problem \eqref{eq5.1} has a solution $u\in W^{1,p}_0 (\mathbb
R^N)$, which is unique up to an element in ${\mathcal P}_{[1-N/p]}$ and
satisfies the estimate
\begin{equation*}
\|u\|_{W^{1,p}_0(\mathbb R^N)/{\mathcal P}_{[1-N/p]}}\le C
\|f\|_{W^{-1,p}_0(\mathbb R^N)}.
\end{equation*}
Moreover, if $1<p< N$, then $u = E*f$. 
\item[\textup{(ii)}]
If $f\in L^p(\mathbb R^N)$, then the problem
\eqref{eq5.1} has a solution $u\in W^{2,p}_0 (\mathbb R^N)$, which is
unique up to an element in ${\mathcal P}_{[2-N/p]}$ and if
$1<p<N/2$, then $ u= E*f $.
\end{enumerate}
\end{theo}

The Calder\'on-Zygmund inequality\index{Calder\'on-Zygmund inequality}
\begin{equation*}
\biggl\|
\frac{\partial^2\varphi}{\partial x_i\partial
x_j}\biggr\|_{L^p(\mathbb R^N)} \le C (N,p) \|\Delta \varphi\|_{L^p(\mathbb
R^N)}, \qquad \varphi \in {\mathcal D}(\mathbb R^N),
\end{equation*}
and Theorem~\ref{Th5.1} are crucial for establishing
Theorem~\ref{Th5.2}. Further, point (i) means that the Riesz
potential of second order satisfies
\begin{equation*}
I_2: W^{-1,p}_0 (\mathbb R^N) \bot {\mathcal P}_{[1-N/p']} \to
W^{1,p}_0 (\mathbb R^N)/{\mathcal P}_{[1-\frac{N}{p'}]}
\end{equation*}
(where the initial space is the orthogonal complement of 
${\mathcal P}_{[1-N/p']}$ in $W^{-1,p}_0 (\mathbb R^N)$)
and it is an isomorphism.

Note that here
\begin{equation*}
W^{1,p}_0 (\mathbb R^N) = \{v\in L^{p^*} (\mathbb R^N); \, \nabla v \in
L^p (\mathbb R^N)\}
\end{equation*}
for $1<p<N$ and ${1}/{p^{*}}={1}/{p}- {1}/{N}$.
And 
for $1< r < \frac{N}{2}$, we have also the continuity property
\begin{equation*}
I_2: L^r (\mathbb R^N) \to L^q (\mathbb R^N), \quad 
   \textup{for $\frac1{q} = \frac1{r} -\frac2{N}$}.
\end{equation*}

\begin{remark}\label{rem5.3}
The problem
\begin{equation}\label{eq5.5}
u -\Delta u = f \qquad \textup{in $\mathbb R^N$}
\end{equation}
is of a completely different nature than the problem \eqref{eq5.1}. The class
of function spaces appropriate for the problem \eqref{eq5.5}
are the classical Sobolev spaces. 
With help of the Calder\'on-Zygmund theory 
one can prove that if $f\in
L^p(\mathbb R^N)$, then the unique solution of \eqref{eq5.5}
 belongs to
$W^{2,p}(\mathbb R^N)$ and can be represented as the Bessel potential
of second order (see \cite{S}): $u = G * f$, where $G$ is the
appropriate Bessel kernel\index{Bessel kernel}, that is, $G$, for
which $\widehat G (\xi) \sim (1+|\xi|^2)^{-1/2}$. Recall that in 
particular $G(x) \sim |x|^{-1} e^{-|x|}$ for $N=3$. 
In the Hilbert case $f\in L^2 (\euc N)$ 
we get
\begin{equation*}
(1+|\xi|^2)\widehat u \in L^2 (\mathbb R^N),
\end{equation*}
which, by Plancherel's theorem, implies that $u\in H^2(\mathbb R^N)$.
For $f\in W^{-1,p} (\mathbb R^N)$, the problem
\ref{eq5.5} has a unique solution $u\in W^{1,p} (\mathbb R^N)$ 
 satisfying the estimate
\begin{equation*}
\|u\|_{W^{1,p}(\mathbb R^N)} \le C (p,n) \|f\|_{W^{-1,p}(\mathbb R^N)}. 
\end{equation*}
\end{remark}

\subsection { Exterior domain}

We consider the problem in an exterior domain with the Dirichlet
boundary condition
\begin{align}\label{eq5.6}\begin{split}
-\Delta u &= f \qquad \textup{in $\Om$},\\
u &= g \qquad \textup{on $\partial\Om$},
\end{split}
\end{align}
where $f\in W^{-1,p}_0(\Om)$ and $g\in
W^{1-1/p,p}(\partial \Om)$.
Invoking the results for $\mathbb R^N$ and bounded domains one
can prove the existence of a solution $u\in W^{1,p}_0 (\Om)$ which
is unique up to an element of the kernel $A^p_0 (\Om) = \{z
\in W^{1,p}_0 (\Om); \, \Delta z = 0\}$ provided that $f$
satisfies the compatibility condition
\begin{equation*}
\langle f,\varphi\rangle 
 = \biggl\langle g, \frac{\partial \varphi}{\partial n}\biggr \rangle
  \qquad \textup{for all $\varphi \in A^{p'}_0 (\Om)$}.
\end{equation*}
The kernel can be characterized in the
following way: It is reduced to $\{0\}$ if $p=2$ or $p<N$ and if
not, then
\begin{equation*}
A_0^p (\Om) = \{C(\lambda-1);\, C\in \mathbb R\} 
\quad \mbox{if} \ p \ge N \ge 3,
\end{equation*}
where $\lambda$ is (unique) solution in $W^{1,2}_0 (\Om) \cap
W^{1,p}_0 (\Om)$ of the problem $\Delta \lambda = 0$ in $\Om$
and $\lambda = 1$ on $\partial \Om $, and
\begin{equation*}
A^p_0 (\Om) = \{C(\mu-u_0); \ C\in\mathbb R\} \quad 
\mbox{if} \ p>N = 2,
\end{equation*}
where $u_0 (x) = (2\pi |\Gamma |)^{-1}\int_{\Gamma }
\log |y-x|\,d\sigma_y$ and $\mu$ is the only solution in
$W^{1,2}_0 (\Om) \cap W^{1,p}_0(\Om)$ of the problem $\Delta \mu =
0$ in $\Om$ and $\mu=u_0$ on $\Gamma $.

\begin{remark}\label{rem5.4}
Similar results exist for the Neumann problem in an exterior
domain (see \cite{AGG}).
The framework of the
spaces  $W_{\alpha}^{m,p}(\mathbb R^N_+)$ also for the Dirichlet
problem in  $\mathbb R^N_+$ was considered in the literature, too.
For a more general theory see \cite{KM}.  
\end{remark}

\section{Elliptic systems}\index{elliptic systems}

\subsection{The Stokes system}\index{Stokes systems}

The Stokes problem is a classical example in the fluid mechanics.
This system models the slow motion with the field of the velocity
$\vec u$ and the pressure $\pi$, satisfying
\begin{equation*}
\begin{array}{rcl}
-\nu \Delta \vec u + \nabla \pi &= &\vec f \qquad\mbox{in} \quad
   \Om,\nonumber\\ 
\operatorname{div} \vec u& =& h \qquad\mbox{in}\quad
   \Om, 
       \tag{S}\\ 
\vec u &= &\vec g \qquad\mbox{on} \quad \Gamma = \partial\Om , \nonumber
\end{array}
\end{equation*}
where $\nu >0$ denotes the viscosity, $\vec f$ is an exterior force,
$\vec g$ is the velocity of the fluid on the domain boundary and $h$ measures 
the
compressibility of the fluids (if $h=0$, it is an incompressible
fluid). The functions $h$ and $g$ must satisfy the compatibility
condition
\begin{equation} \label{eq4.1}
\int_{\Om} h(x)\,dx = \int_{\Gamma } \vec g \cdot
\vec n \, d\sigma. 
\end{equation}

\begin{theo}\label{Th4.1}
Let $\Om$ be a Lipschitz bounded domain in
$\mathbb R^N$, $N \ge 2$. Let $\vec f \in H^{-1} (\Om)^N$, $h\in
L^2(\Om)$ and $g \in H^{1/2}(\Gamma)^N$ satisfy \eqref{eq4.1}. Then the
problem $(S)$ has a unique solution $(\vec u, \pi) \in
H^1(\Om)^N\times L^2(\Om)/\mathbb R$ satisfying the a priori 
estimate
\begin{equation*}
\|\vec u\|_{H^1(\Om)} + \| \pi\|_{L^2(\Om)/\mathbb R} \le C(\|\vec
f\|_{H^{-1}(\Om)} + \|h\|_{L^2(\Om)}
+\|g\|_{H^{1/2}(\Gamma)}).
\end{equation*}
\end{theo}

\noindent 
In order to prove Theorem~\ref{Th4.1} one can start with a homogeneous
problem. The procedure of finding $\vec u$ is a simple application of
the Lax-Milgram theorem. 
Application of De Rham's theorem gives the pressure $\pi$. We
introduce the space 
\begin{equation*}
{\mathcal V} 
 =\{\vec v \in{\mathcal D}(\Om)^N;\,\operatorname{div}\vec v=0\}
\end{equation*}
and define $\vec F \in H^{-1} (\Om)^N$ by
\begin{equation*}
 \langle \vec F, \vec v\rangle_{H^{-1}\times \oH^1} = 0 \quad 
   \textup{for all  $\vec v \in {\mathcal V}$}.
\end{equation*}
Moreover, there exists $\pi \in L^2(\Om)$, unique up to an
additive constant, and such that $\vec F = \nabla \pi$. The problem
$(S)$, which we transform to the homogeneous case ($h=0$, $g=0$), can be
formulated on an abstract level. Let $X$ and $M$ be two real Hilbert
spaces and consider the following variational problem:
Given $\vec L \in X'$ and $X\in M'$, find $(\vec u,\pi) \in X \times
M$ such that
\begin{align}\label{eq4.2}\begin{split}
{\mathcal A} (\vec u, \vec v) + B [\vec v, \pi] &= \vec L (\vec v), 
   \qquad \ v\in X,\\
B[\vec u,q] &=X(q),  \qquad  q\in M.
\end{split}
\end{align}
where the bilinear forms ${\mathcal A}$, $B$ and the linear form $\vec
L$ are defined by
\begin{equation*}
{\mathcal A}(\vec u,\vec v)= \int_{\Om }\nabla \vec u \cdot
\nabla \vec v,
\quad B[\vec v,q]=-\int_{\Om }[q\nabla \cdot
 \vec v] , 
\quad L(\vec v)=\int_{\Om }\vec f\cdot \vec v .
\end{equation*}

\begin{theo}\label{Th4.3}
If the bilinear form ${\mathcal A}$ is
coercive\index{coercivity} in the space 
\begin{equation*}
V = \{\vec v \in X; \, B[\vec v, q]=0\ \textup{for all $q \in M$},
\end{equation*}
i.e.\ if there exists $\alpha>0$ such that
\begin{equation*}
{\mathcal A} (\vec v,\vec v)
  \ge \alpha\|\vec v\|^2_X, \qquad \vec v \in V,
\end{equation*}
then the problem \eqref{eq4.2} has a unique solution $(\vec u,\pi)$ if and
only if the bilinear form ${\mathcal B}$ satisfies the ``\,$\inf$-$\sup$''
condition:\index{inf-sup condition}
\begin{equation*}
\mbox { there exists }\beta > 0 \mbox { such that }\ \inf_{q\in
M} \sup_{\vec v \in X} \frac{{\mathcal B}(\vec v, 
q)}{\|\vec v\|_X \|q\|_M} \ge \beta.
\end{equation*}
\end{theo}

\noindent As for the Dirichlet problem, the regularity result is the
following:\index{regularity}

\begin{theo}\label{Th4.3a}
Let $\Om$ be a bounded domain in $\mathbb R^N$, of
the class $C^{m+1,1}$ if $m\in\mathbb N$ and $C^{1,1}$ if $m=-1$. Let
$f\in W^{m,p}(\Om)^N$, $h\in W^{m+1,p} (\Om)$ and $\vec g \in
W^{m+2-1/p,p}(\Gamma)^N$  satisfy condition
\eqref{eq4.1}. Then the problem $(S)$ has a unique solution
$(\vec u,\pi) \in W^{m+2,p}(\Om)^N\times W^{m+1,p}(\Om)/\mathbb R$.
\end{theo}

\begin{remark}\label{rem4.4}
It is possible to solve $(S)$ under weaker assumption, for instance
if $\vec f \in W^{-1/p}(\Om')$, $h=0$ and $\vec g \in
W^{-1/{p},p} (\Gamma)^N$. We can prove that then $(\vec u,\pi)
\in L^p (\Om)^N\times W^{-1,p}(\Om)$.
\end{remark}

\subsection{The linearized elasticity}\index{linearized elasticity}
The equations governing
the displacement $ \vec u = (u_1,u_2,u_3)$ of a three
dimensional structure subjected to an external force field
$\vec{f}$ are written as ($\Om$ is a bounded open subset
of $\euc 3$ and $\Gamma = \partial
\Om$)
\begin{align*}
-\mu  \Delta \vec u -(\lambda+\mu ) 
   \nabla (\nabla \cdot \vec u) &= \vec f\qquad \hbox{ in }\Om,\\
\vec u &= 0  \qquad \hbox{ on }\Gamma_0,  \\
\sum_{j=1}^3 \sigma_{ij}(\vec u ) \nu_j &=\vec g_i  
   \quad \quad \hbox{on }\Gamma_1=\Gamma-\Gamma_0,  
\end{align*}
where
$\lambda >0$ and $\mu >0$ are two material
characteristic constants, called the \textit{Lam\'e coefficients}, and
($ \vec v  = (v_1,v_2,v_3)$)
\begin{align}\label{DefSigma}\begin{split}
\sigma_{ij}(\vec v)  = \sigma_{ji}(\vec v) & = \lambda
 \delta_{ij} \, \sum_{k=1}^3\varepsilon_{kk}(\vec v)
  +2\mu\varepsilon_{ij}(\vec v)\\
\textup{with\ } \varepsilon_{ij}(\vec v) &= \varepsilon_{ji}(\vec v) 
  =\frac{1}{ 2} (\partial_j v_i + \partial_i v_j),
\end{split}
\end{align}
where $\delta_{ij}$ denotes the Kronecker symbol, i.e.\  $\delta_{ij}
=1$, for $i=j$ and  $\delta_{ij} 
=0$, for $i \ne j$. These equations describe the equilibrium of
an elastic homogeneous isotropic body that cannot move along
$\Gamma_0$; along $\Gamma_1$, surface forces of density $\vec g
=(g_1,g_2,g_3)$ are given. The case $\Gamma_1 =
\emptyset$ physically corresponds  to clamped structures. The
matrix with entries $\varepsilon_{ij}(\vec u)$ is the linearized
strain tensor while $\sigma_{ij} (\vec u )$ represents the
linearized \index{linearized stress tensor} stress tensor; the relationship 
\eqref{DefSigma}
between these tensors is known as \textit{Hooke's law}.%
\index{Hooke's law}
We refer for example to \cite{CiarletLions}, 
\cite{NH} (and to the references herein) for most of the results
stated in this paragraph. The variational formulation of this problem is
\begin{equation}\label{PVElas}
\textit{to find $\vec u \in V$ such that  ${\mathcal
A}(\vec u,\vec v)=L(\vec v)$ for all $\vec v\in V$},
\end{equation}
where the bilinear form ${\mathcal A}$ and the linear form $L$ are given
by
\begin{subequations}\label{DefALElas}
\begin{equation}\label{DefAElas}
{\mathcal A}(\vec u,\vec v) = \int_\Om [ \lambda (\nabla
 \cdot\vec u) (\nabla \cdot \vec v)
  +2\mu\sum_{i,j=1}^3 \varepsilon_{ij} (\vec u) \varepsilon_{ij}
   (\vec v)](x)\,dx,
\end{equation}
\begin{equation}\label{DefLElas}
L(\vec v)=\int_\Om \vec f(x) \cdot \vec v(x)\,dx
  +\int_{\Gamma_1} \vec g(x) \cdot \vec v(\sigma)\, d\sigma;
\end{equation}
\end{subequations}
\bigskip
The functional space $V$ is defined as
\begin{equation*}
V=\{\vec v =(v_1,v_2,v_3) \in [H^1(\Om)]^3; 
 \, \gamma_0 v_i = 0 \,\hbox{ on }\Gamma_0, \ 1\le i\le3\}.
\end{equation*}
To prove the ellipticity of ${\mathcal A}$, one needs the following
Korn inequality: There exists a positive constant $C(\Om)$
such that, for all $\vec v =(v_1,v_2,v_3)
\in[H^1(\Om)]^3$, we
have
\index{Korn's inequality}
\begin{equation}\label{IneKorn}
\| \vec v |
|_{1,\Om}\le C(\Om) \, \left
[\sum_{i,j=1}^3 \| \varepsilon_{ij}(\vec v)\|
 _{L^2(\Om)}^2+\sum_{i=1}^3 \| v_i
  \| _{L^2(\Om)}^2\right]^{1/2}.
\end{equation} 

The following result holds true:

\begin{theo}
Let $\Om$ be a bounded open in $\euc 3$ with a
Lip\-schitz boundary, and let $\Gamma_0$ be a measurable subset of
$\Gamma$, whose measure (with respect to the surface measure $d
\Gamma(x)$) is positive. Then the mapping
\begin{equation*}
\vec v \mapsto \left [\sum_{i,j=1}^3
\| \varepsilon_{ij}(\vec v)\|_{L^2(\Om)}^2\right ]^{1/2}
\end{equation*}
is a norm on $V$, equivalent to the usual norm 
$\| \,.\, \|_{1,\Om}$.
\end{theo}

As a consequence, we get:

\begin{theo}\label{ThmExisElas}
Under the above assumptions, there exists a unique $u \in V$ solving
the variational problem \eqref{PVElas}-\eqref{DefALElas}. This
solution is also the unique one which minimizes 
the energy functional 
\begin{align*}
E(\vec v) &=  \frac{1 }{ 2} \,
\int_\Om [\lambda (\nabla \cdot \vec v)^2  + 
2 \mu\sum_{i,j=1}^3 [\varepsilon_{ij} (\vec v)]^2](x)\,dx \\
&\quad
-\left[\,\int_\Om \vec f(x) \cdot \vec v(x)\,dx + 
\int_{\Gamma_1}\vec g(x) \cdot \vec v(\sigma) \,d\sigma \, \right
]
\end{align*}
over the space $V$.
\end{theo}

\noindent 
{\bf Acknowledgment:} M.~Krbec and \v S.~Ne\v casov\' a were
supported by the Institutional Research Plan No.\ AV0Z10190503 and
by the Grant Agency of the Academy of Sciences No.\ IAA10019505.

\printindex
\end{document}